# Similarity Analytical Solutions for the Schrödinger Equation with the Riesz Fractional Derivative in Quantum Mechanics

**Asim Patra**

*National Institute of Technology*
*Department of Mathematics*
*Rourkela-769008, India*
*Email: asimp1993@gmail.com*

**Abstract**

The present article deals with the similarity method to tackle the fractional Schrödinger equation where the derivative is defined in the Riesz sense. Moreover the procedure of reducing a fractional partial differential equation (FPDE) into an ordinary differential equation (ODE) has been efficiently displayed by means of suitable scaled transform to the proposed fractional equation. Furthermore the ODEs are treated effectively via the Fourier transform. The graphical solutions are also depicted for different fractional derivatives $\alpha$.

## 1. Introduction

The concept of fractional calculus has aroused interests among the researchers since a very long time and emerged as a desirable field on account of its applications in numerous areas like engineering sciences, social sciences, physics, mathematics and even financial matters and modelling for the real life processes[1-4]. An enormous number of proficient ideas have been emerging in the particular field from the beginning of its inception. The most significant advantage of the fractional order models is that they can give a fluent explanation regarding the hereditary characteristics of various processes such as anomalous diffusion, viscoelasticity, and many other [5]. Further the concept of fractional derivatives includes several types of derivatives like the Riemann-Lioville(R-L) derivative, Riesz derivative, Caputo derivative and other which has numerous applications in various areas in the fractional sense. The concept of Riesz fractional derivative is one of the crucial elements of the space fractional quantum mechanics that has been studied in the infinite as well as in the finite domain by several researchers which is mathematically defined as the amalgamation of the left R-L derivative and right R-L derivative.

For the purpose of investigation of the analytical solutions of the FPDEs, they are basically solved by the separation of variable method [1] or by employing the Laplace transform along with the Fourier transformation. Further certain other semi-analytic techniques such as the homotopy perturbation method[6], homotopy analysis method [7], variational iteration method [8], fractional differential transformation method[9,10], Adomian decomposition method[11,12] etc. have been recently employed to treat the FPDEs for obtaining their series solutions. Apart from this another analytic technique viz. Lie group method is utilized in [13] for solving the space-time fractional diffusion equation.

In the context of the present paper, the Riesz space fractional Schrödinger equation with variable coefficients given by

$$\varepsilon \frac{\partial u}{\partial t} - i \frac{\varepsilon^2}{2} t^{\alpha/2} \frac{\partial^\alpha u}{\partial |x|^\alpha} = -iV(x)u(x,t) \quad (1.1)$$

has been considered for obtaining the similarity solutions where $\frac{\partial^\alpha u}{\partial |x|^\alpha}$ is the Riesz derivative of the auxiliary wave function $u(x,t)$ for fractional order $\alpha$. Further $V(x)$ denotes the electrostatic potential function which is taken as unity in our case. The fractional Schrödinger equation is employed in the quantum field theory and quantum mechanics for the purpose of evolution of wave packets. There are many researchers who had tackled the equation via numerical approach. So it will be a worthy and interesting task to handle the fractional version of the Schrödinger equation via a novel analytical approach.

In the present research article, an efficient analytical technique is utilized to solve the Riesz fractional Schrödinger equation (RFSE) involving the Riesz fractional derivative. The similarity method is highly effective while treating the linear as well as non-linear partial differential equations (PDEs). In the literature survey certain numerical techniques are provided for displaying the numerical simulation for the results obtained for Schrödinger equation but they are not as much accurate as the results obtained by an analytical method. By deriving the motivation from the accuracy of the analytical solutions, the similarity solutions are calculated and manifested for the RFSE. In other words, a direct technique is employed here in which first of all the fractional PDE is transformed to a fractional ODE by suitable scaling transformation. Then the resulting ODE is solved by means of Fourier transformation. Further graphical solutions are depicted for an effective comparative study between the solutions of fractional order and integer order.

The literature regarding the analytical treatment of fractional Schrodinger equation is limited to certain extent. For this reason the analytical technique illustrated for the proposed problem in this article would prove to be highly efficient. Some of the researchers have tackled the integer order Schrodinger equation via certain numerical approach. Elsaid et al[14] have illustrated the similarity solutions for PDEs involving the Riesz derivative. Furthermore Elsaid et al[15], in a different paper, have tackled with the fractional heat equation in context with the similarity solutions. Djordjevic et al[16] have utilized the similarity method for the nonlinear fractional Korteweg-de Vries equation. Apart from that the symmetry properties has been discussed for the fractional diffusion equation by Gazizov et al[17]. Further the scaling transformation have been studied by Kandasamy et al[18] for the fluid viscosity on the process of mass and heat transfer. Apart from that El-Sheikh[19] have studied regarding the similarity solutions for KdV equation via symmetry method. In a different paper, the Burgers equation had been treated in terms of freezing solutions by Matthes [20].

The paper is organized in the prescribed way: First in sec 2 some definitions relating to the Riesz fractional derivative has been manifested which is required in the later section. Further the similarity technique has been elaborately displayed for the proposed problem in the next section. In other words, this section deals with the application of the similarity method for the space fractional Schrödinger equation by means of the scaled transformation thus converting the FPDE to fractional ODE. Further in sec 4, the fractional ODE is solved by utilizing the Fourier transform technique. In the next section certain graphical solutions are depicted for the resulting solution of the fractional equation for different values of fractional order $\alpha$ and compared with the integer order. Furthermore the conclusion is displayed in the final section 5.

## 2. Basic preliminaries of fractional calculus

In this section, the basic concepts regarding the Riesz derivative have been discussed. Moreover some of the formulae, which will be later used, regarding the discrete Fourier transform and inverse discrete Fourier transform has been displayed.

**Definition I(a)**: The Riesz derivative or Riesz-Feller derivative of a function $f(x)$ of fractional order $\alpha$, denoted by $R_x^\alpha f(x)$, is usually defined in terms of its Fourier transform as

$$R_x^\alpha f(x) = \mathcal{F}^{-1}[-|k|^\alpha \hat{f}(k)], \tag{2.1}$$

where $|k|^\alpha$ is the Fourier symbol for $k \in R$ and $\hat{f}(k)$ is the Fourier transform of $f(x)$ given by

$$\hat{f}(k) = \mathcal{F}(f(x)) = \int_{-\infty}^{\infty} e^{ikx} f(x) dx. \tag{2.2}$$

**Definition I(b)**: The Riesz derivative can also be defined in other way. For a fractional order $\alpha$ $(n-1 < \alpha \leq n)$, the Riesz derivative on the infinite domain $-\infty < x < \infty$ is defined as [9]

$$-(-\Delta)^{\alpha/2} \psi(x) = \frac{\partial^\alpha \psi(x)}{\partial |x|^\alpha} = -c_\alpha \left( {}_{-\infty}D_x^\alpha \psi(x) + {}_xD_\infty^\alpha \psi(x) \right), \tag{2.3}$$

where ${}_{-\infty}D_x^\alpha \psi(x)$ is the Riemann-Liouville left derivative defined as

$${}_{-\infty}D_x^\alpha \psi(x) = \frac{1}{\Gamma(n-\alpha)} \frac{\partial^n}{\partial x^n} \int_{-\infty}^{x} \frac{\psi(\xi)d\xi}{(x-\xi)^{1-n+\alpha}}, \tag{2.4}$$

and ${}_xD_\infty^\alpha \psi(x)$ is the Riemann-Liouville right derivative defined as

$${}_xD_\infty^\alpha \psi(x) = \frac{(-1)^n}{\Gamma(n-\alpha)} \frac{\partial^n}{\partial x^n} \int_{x}^{\infty} \frac{\psi(\xi)d\xi}{(\xi-x)^{1-n+\alpha}}, \tag{2.5}$$

$$c_\alpha = \frac{1}{2\cos\frac{\alpha\pi}{2}}, \quad \alpha \neq 1. \tag{2.6}$$

There is a special case when $\alpha = 1$, in which the derivative is taken for the Hilbert transform $H$ of the function $\psi(x)$ and the integral is considered in the sense of Cauchy principal value given by

$$D_x^1 \psi(x) = \frac{d}{dx} H\psi(x) = \frac{d}{dx} \frac{1}{\pi} \int_{-\infty}^{\infty} \frac{\psi(z)dz}{z-x}, \tag{2.7}$$

If $\psi(x,t)$ is given on the finite domain $[a,b]$ and it satisfies the boundary conditions $u(a,t)=u(b,t)=0$, the function can be extended by taking, for $x \leq a$ and $x \geq b$, $\psi(x,t) \equiv 0$. Thus, the Riesz fractional derivative of order $\alpha$ $(n-1<\alpha \leq n)$ on the finite interval $a \leq x \leq b$ can be defined as

$$\frac{\partial^\alpha \psi(x,t)}{\partial |x|^\alpha} = -\frac{1}{2\cos\frac{\alpha\pi}{2}}\left({}_a D_x^\alpha \psi(x,t) + {}_x D_b^\alpha \psi(x,t)\right), \tag{2.8}$$

where

$$_a D_x^\alpha \psi(x,t) = \frac{1}{\Gamma(n-\alpha)}\frac{\partial^n}{\partial x^n}\int_a^x \frac{\psi(\xi,t)d\xi}{(x-\xi)^{1-n+\alpha}}, \tag{2.9}$$

$$_x D_b^\alpha \psi(x,t) = \frac{(-1)^n}{\Gamma(n-\alpha)}\frac{\partial^n}{\partial x^n}\int_x^b \frac{\psi(\xi,t)d\xi}{(\xi-x)^{1-n+\alpha}}. \tag{2.10}$$

## 3. Similarity method for the Schrödinger equation involving the Riesz fractional derivative

The section involves the analytical treatment via the similarity method solutions to the following fractional Schrödinger equation

$$\varepsilon t\frac{\partial u(x,t)}{\partial t} - i\frac{\varepsilon^2}{2}t^{\alpha/2}\frac{\partial^\alpha u(x,t)}{\partial |x|^\alpha} = -iV(x)u(x,t), \tag{3.1}$$

where $\frac{\partial^\alpha u(x,t)}{\partial |x|^\alpha}$ denotes the Riesz definition of fractional derivative of order $\alpha$ which can also be illustrated as $-(-\Delta)^{\alpha/2} u(x,t) = \frac{\partial^\alpha u(x,t)}{\partial |x|^\alpha}$.

Let us take the variables according to the scaling transformations in the following way,

$$t = \mu^n \tilde{t}, \quad x = \mu^p \tilde{x}, \quad u = \mu^q \tilde{u}, \tag{3.2}$$

where $\mu$ is the scaling factor.

The first term of eq. (3.1) is transformed as

$$\varepsilon t\frac{\partial u}{\partial t} = \varepsilon\mu^n \tilde{t}\frac{\partial(\mu^q \tilde{u})}{\partial \tilde{t}}\frac{\partial \tilde{t}}{\partial t} = \varepsilon\mu^n \tilde{t}\mu^{q-n}\frac{\partial \tilde{u}}{\partial \tilde{t}}. \tag{3.3}$$

By the definition of the Riesz derivative of fractional order $\alpha$ $(n-1<\alpha<n)$, we have

$$\frac{\partial^\alpha u(x,t)}{\partial |x|^\alpha} = -\frac{1}{2\cos\left(\frac{\alpha\pi}{2}\right)}\left({}_{-\infty} D_x^\alpha u(x,t) + {}_x D_\infty^\alpha u(x,t)\right), \quad \alpha \neq 1. \tag{3.4}$$

For $n=2$, we have the Riemann-Liouville left fractional derivative as

$$_{-\infty} D_x^\alpha u(x,t) = \frac{1}{\Gamma(2-\alpha)}\frac{\partial^2}{\partial x^2}\int_{-\infty}^x (x-\xi)^{1-\alpha}u(\xi)d\xi, \tag{3.5}$$

which undergoes transformation (3.2) and further taking $\xi = \mu^p \tilde{\xi}$ we get,

$$_{-\infty}D_{\tilde{x}}^\alpha u(\tilde{x},\tilde{t}) = \frac{\mu^{-2p}}{\Gamma(2-\alpha)}\frac{d^2}{dx^2}\int_{-\infty}^{\tilde{x}}\left(\mu^p\tilde{x} - \mu^p\tilde{\xi}\right)^{1-\alpha}\mu^q\mu^p\tilde{u}d\tilde{\xi}. \tag{3.6}$$

Similarly the Riemann-Liouville right derivative for $n = 2$ is given by

$$_xD_\infty^\alpha u(x,t) = \frac{1}{\Gamma(2-\alpha)}\frac{\partial^2}{\partial x^2}\int_{-\infty}^{x}(\xi - x)^{1-\alpha}u(\xi)d\xi, \tag{3.7}$$

which renders us

$$_{\tilde{x}}D_\infty^\alpha u(\tilde{x},\tilde{t}) = \frac{\mu^{-2p}}{\Gamma(2-\alpha)}\frac{d^2}{dx^2}\int_{\tilde{x}}^{\infty}\left(\mu^p\tilde{\xi} - \mu^p\tilde{x}\right)^{1-\alpha}\mu^q\mu^p\tilde{u}d\tilde{\xi}, \tag{3.8}$$

after the same transformation (3.2).

So the Riesz definition of fractional derivative is transformed as

$$\frac{\partial^\alpha u(x,t)}{\partial |x|^\alpha} = \frac{-1}{2\cos\left(\frac{\alpha\pi}{2}\right)}\frac{\mu^{-2p}}{\Gamma(2-\alpha)}\left(\frac{d^2}{dx^2}\int_{-\infty}^{\tilde{x}}\left(\mu^p\tilde{x} - \mu^p\tilde{\xi}\right)^{1-\alpha}\mu^q\mu^p\tilde{u}d\tilde{\xi} + \frac{d^2}{dx^2}\int_{\tilde{x}}^{\infty}\left(\mu^p\tilde{\xi} - \mu^p\tilde{x}\right)^{1-\alpha}\mu^q\mu^p\tilde{u}d\tilde{\xi}\right)$$

$$= \mu^{-\alpha p + q}\frac{\partial^\alpha \tilde{u}(\tilde{x},\tilde{t})}{\partial |\tilde{x}|^\alpha}. \tag{3.9}$$

Thus substituting the equations (3.3) and (3.9) into (3.1), we obtain

$$\varepsilon\mu^n\tilde{t}\mu^{q-n}\frac{\partial\tilde{u}}{\partial\tilde{t}} - i\frac{\varepsilon^2}{2}\mu^{n\alpha/2}\mu^{-\alpha p+q}\frac{\partial^\alpha \tilde{u}(\tilde{x},\tilde{t})}{\partial|\tilde{x}|^\alpha} = -i\mu^q\tilde{u}. \tag{3.10}$$

Subsequently it is evident from eq. (3.10) that

$$n + q - n = \frac{n\alpha}{2} - \alpha p + q = q, \tag{3.11}$$

from which specifying the condition

$$p = \frac{n}{2}, \tag{3.12}$$

eq. (3.1) is invariant under the scaling transformation (3.2). Further the characteristic equation related to the transformation (3.2) is represented by

$$\frac{dt}{tn} = \frac{dx}{xp} = \frac{du}{uq}. \tag{3.13}$$

Then eq. (3.13) has been solved by taking $q = 0$ which gives the solution expressible in the form

$$u(x,t) = c_1 \text{ and } c_2 = \frac{x}{\sqrt{t}},$$

which provides us with
$$u(x,t) = \varphi(xt^{-1/2}) = \varphi(v) \tag{3.14}$$

where $v = \dfrac{x}{t^{\frac{p}{n}}}$ from eq. (3.12).

Now accordingly converting eq. (3.10) by applying the solution (3.14), we get
$$\varepsilon t \frac{\partial u}{\partial t} = \varepsilon t \frac{\partial(\varphi(v))}{\partial v}\frac{\partial v}{\partial t} = \varepsilon\left(-\frac{p}{n}\right)xt^{-\frac{p}{n}}\frac{\partial \varphi}{\partial v} \tag{3.15}$$

Further transforming eq. (3.9) in terms of the single variable function $\varphi(v)$, we have

$$\frac{\partial^\alpha u(x,t)}{\partial |x|^\alpha} = \frac{-1}{2\cos\left(\frac{\alpha\pi}{2}\right)}\frac{1}{\Gamma(2-\alpha)}\left(t^{-\frac{2p}{n}}\frac{d^2}{dv^2}\int_{-\infty}^{v}\left(vt^{\frac{p}{n}} - yt^{\frac{p}{n}}\right)^{1-\alpha}\varphi(v)t^{\frac{p}{n}}dy + t^{-\frac{2p}{n}}\frac{d^2}{dv^2}\int_{v}^{\infty}\left(yt^{\frac{p}{n}} - vt^{\frac{p}{n}}\right)^{1-\alpha}\varphi(v)t^{\frac{p}{n}}dy\right)$$

$$= t^{-\frac{p}{n}}\frac{\partial^\alpha \varphi(v)}{\partial |v|^\alpha}. \tag{3.16}$$

where $\dfrac{d^2}{dx^2} = t^{-\frac{2p}{n}}\dfrac{d^2}{dv^2}$ and $y = \xi t^{-\frac{p}{n}}$.

Then substituting all the terms in eq. (3.1), we obtain the resulting FDE as
$$-\frac{\varepsilon}{2}v\frac{d\varphi}{dv} - i\frac{\varepsilon^2}{2}\frac{\partial^\alpha \varphi(v)}{\partial |v|^\alpha} = -i\varphi(v). \tag{3.17}$$

Now taking the Fourier transformation both sides of the similar equation as (3.17) viz.
$$-\frac{\varepsilon}{2}t\frac{d\psi}{dt} - i\frac{\varepsilon^2}{2}\frac{\partial^\alpha \psi(t)}{\partial |t|^\alpha} = -i\psi(t) \tag{3.18}$$

we get,
$$\frac{\varepsilon}{2}\left(\omega\frac{d\Psi}{d\omega} + \Psi(\omega)\right) + i\frac{\varepsilon^2}{2}|\omega|^\alpha \Psi(\omega) = -i\Psi(\omega). \tag{3.19}$$

Further solving eq. (3.19) we obtain the solution in terms of Fourier transform as
$$\Psi(\omega) = c\omega^{\frac{-2i}{\varepsilon}-1}e^{-i\varepsilon\frac{|\omega|^\alpha}{\alpha}} \tag{3.20}$$

In the next step, we apply the inverse Fourier transformation for eq. (3.20) and further taking $c = -1/i$, the solution is obtained as
$$\psi(t) = \frac{1}{2\pi}\int_{0}^{t}\int_{-\infty}^{\infty}e^{-i\frac{\varepsilon|\omega|^\alpha}{\alpha}}e^{-i\omega\frac{-2i-\varepsilon}{\varepsilon}\tau}d\omega d\tau \tag{3.21}$$

which can be also denoted, taking $\varepsilon = 2$, as

$$\phi(v) = \frac{1}{2\pi}\int_0^v \int_{-\infty}^{\infty} e^{-i\frac{|\omega|^\alpha}{\alpha}} e^{-i\omega^{-i-1}\tau} d\omega\, d\tau. \tag{3.22}$$

After a computer calculation via a reliable software mathematica, we get different resulting functions for various values of $\varepsilon$ and order $\alpha$ in terms of standard error function.

We get the plot for different values of $\alpha$ for $\alpha = 1.7, \alpha = 1.9$ and the classical order $\alpha = 2.0$.

## 4. Graphical solutions and discussions

In this section the graphical solutions are displayed for the function $\phi(v)$ for different values of fractional order $\alpha$. Further the effect of changing the fractional order has been efficiently depicted on the behaviour of the resulting function $\phi(v)$ obtained by solving the fractional Schrodinger equation. Apart from that it has been shown that when the order $\alpha$ approach to 2, the solution depicted graphically for the equation(3.1) displays the error function of classical order. Here in Fig 1 the graph has been displayed for similarity method solutions for eq. (3.1) using different values of $\alpha$ and $\varepsilon = 1.2$. In Fig 2, Fig 3 and Fig 4, the graphical solutions has been depicted for $\varepsilon = 2$, $\varepsilon = 3$ and $\varepsilon = 10$ respectively.

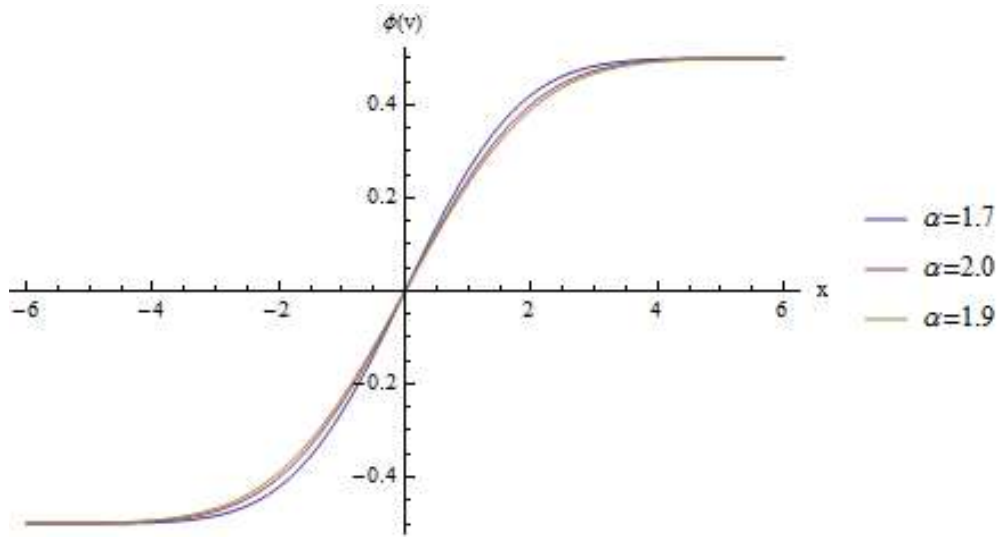

**Fig 1.** Graphical solution for eq. (3.1) obtained by similarity solutions method at different values of $\alpha$ taking $\varepsilon = 1.2$.

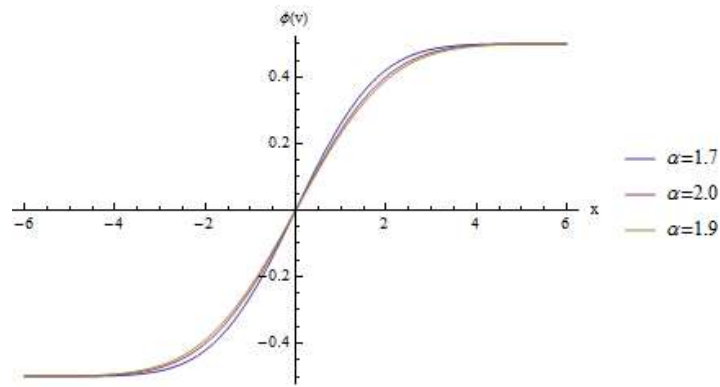

**Fig 2.** Graphical solution for eq. (3.1) obtained by similarity solutions method at different values of $\alpha$ taking $\varepsilon = 2$.

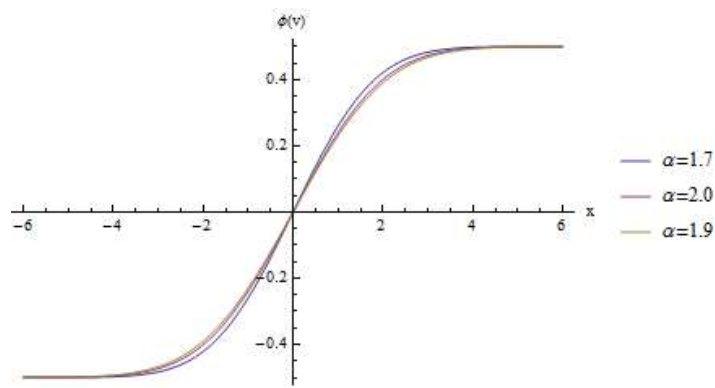

**Fig 3.** Graphical solution for eq. (3.1) obtained by similarity solutions method at different values of $\alpha$ taking $\varepsilon = 3$.

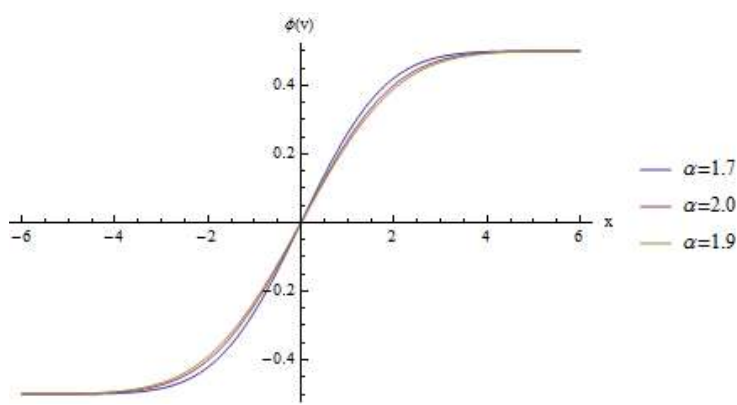

**Fig 4.** Graphical solution for eq. (3.1) obtained by similarity solutions method at different values of $\alpha$ taking $\varepsilon = 10$.

**Conclusion**

The paper deals with the similarity solution method for treating the fractional Schrodinger equation where the fractional derivative is taken in the Riesz sense. The proposed technique transforms the fractional equation with two or more independent variables to single variable differential equation form and thereby easing the task of finding the required solutions. The transformed differential equation is then solved by utilizing the Fourier transformation method. Apart from this, the resulting solution function has been graphically depicted for several values of fractional order $\alpha$ and comparing with the integer order $\alpha$ when $\alpha$ tends to 2.

**Acknowledgement**

The authors take the opportunity to express his heartfelt gratitude to the learned anonymous reviewer for his valuable comments and suggestions for the improvement and betterment of the manuscript.

**Conflict of interest**

The author declares that there is no conflict of interest regarding the publication of this research paper.